
\input vanilla.sty
\magnification 1200
\overfullrule=0pt
\baselineskip 18pt
\input definiti.tex
\input mathchar.tex
\define\pmf{\par\medpagebreak\flushpar}

\define\e{\varepsilon}
\define\vp{\varphi}
\define\pbf{\par\bigpagebreak\flushpar}

\define\k{\kappa}

\font\frak=eufm10
\pmf
\title
Volumes of Discrete Groups and Topological Complexity of Homology Spheres
\endtitle
\author
Alexander Reznikov
\endauthor

Can one ascribe a ``size'' to a discrete group?  By that we mean
a real-valued invariant, $\mu (\Gamma)$ satisfying the following condition: for
any subgroup $\Delta < \Gamma$ of finite index,
$\mu (\Delta) = (\Gamma : \Delta) \mu (\Gamma)$.  The importance of
such an invariant for
the structure of $\Gamma$ is clear from the following observation: if
$\mu (\Gamma) > 0$ and
$\Delta_1, \Delta_2$ are two isomorphic subgroups of
$\Gamma$ of finite index,
then $(\Gamma: \Delta_1) = (\Gamma: \Delta_2)$.

Some rational-valued invariants satisfying the axiom above, have been
studied by Wall [17], Bass [1], Chiswell [7], Brown [4], Stallings [16] etc.
and called (generalized) Euler characteristics.  These are defined
for a FP group
$\Gamma$ by $\chi (\Gamma) = \sum (-1)^i r k P_i$ where $\bbp \to \bbq$ is a
finite f.g. projective
resolution of $\bbq$ and $rk: K_0 (\bbq \Gamma) \to \bbz$ is a
homomorphism and extended to all VFP groups in the usual way [17].

On the other hand, if $\Gamma$ runs through fundamental groups of
closed three-manifolds, there are two known invariants, (usually called
covering invariants by obvious reason).  One is Gromov's simplicial
volume, which is positive for any hyperbolic manifold.  The other is
Wang and Wu's invariant [18], defined for any plumbing of Seifert fibrations.

The approach presented in this paper, takes a new standpoint in
dealing the above problem,
that is, a ``measure-theoretic'' one.  Namely, we start
with ``semimultiplicative invariants'', analogous to (outer)
measure in measure theory and then produce a truly
multiplicative invariants by
going to lower limit over the directed set of finite index subgroups.

In course of application, we address two fundamental and
well-known problems of Gromov and Lyndon:

\demo{Problem A} (Gromov, see [5]).  Consider a category $M_n$ of
closed manifolds of
dimension $n$ with nonzero-degree ways as morphisms.  Study a partial
order $M \ge N \Leftrightarrow Mor (M, N) \neq \phi$.  For which $N$ the
degrees of maps $f: M \to N$ are bounded for all $M$?

\demo{Problem B} (Lyndon, [12], problem 13).  Extend and relate the
theories of
deficiency, the rate of growth and the Euler-Poincar\'e characteristic.  In
particular, what influence does the deficiency
have on the structure of an infinite group?

As an answer to the second problem, we will show how to
construct the \underbar{deficiency volume}
$V_d (G)$, which will imply the following

\proclaim{Theorem 1}  Let $G$ be a f.p. group with def$(G)
\ge 2$.  Then any two
isomorphic finite-index subgroups of $G$ have the same index.
\endproclaim

We will then introduce the rank volume $V_r (G)$ and
will prove that $V_r (G) > 0$ if
def$(G) \ge 2$.  We will show that if $\Gamma$ maps onto $\Pi$, then
$V_r (\Gamma) \ge V_r (\Pi)$, and derive the following corollary.

\proclaim{Corollary}  Let $M$ be a compact manifold with residually finite
fundamental group of deficiency $\ge 2$.  Then any nonzero degree map $f: M \to
M$ induces
an automorphism of $\pi_1 (M)$.
\endproclaim

We introduce then our second main tool, the volume of a
representation, as presented in [13]
and [14].  For $\Gamma = \pi_1 (M)$ where $M$ is a closed
$n$-dimensional manifold, and a representation
$\rho$ of $\Gamma$ in a real semi-simple Lie group $G$, Vol$(G) =
(\rho^\ast (Bor_n), [M])$, where
$[M] \in H_n (\Gamma, \bbz)$ is an image of the
fundamental class of $M$ and
Bor$_n \in H^n_{\text{cont}} (G)$ is a Borel
generator for the continuous cohomology.
For example, if $G = SL_2 (\bbc)$ then
Vol$(\rho)$ is the hyperbolic volume of $\rho$
[9], [13], [14], and if $G = \widetilde{SL_2 (\bbr)}$, then
Vol$(\rho)$ is the Seifert volume or a
Chern-Simons invariant [2,3].  We set
$\overline{V_G} (\Gamma) = \max\limits_{\rho} \text{Vol} (\rho)$;
this is an lower volume.
Observe that the maximum is finite by Cheeger-Simons
rigidity [6], [14].  If
$G = S O (n, 1)$ and $\Gamma = \pi_1 (M)$, then using ideas of Gromov
and Thurston, one shows [14] $\overline{V}_G (\Gamma) \le \mu_n \cdot
\| M \|$; it follows that
$\overline{V}_G (\Gamma)$ gives arise to a volume in that case.

Let $\Gamma$ be a uniform lattice in $G$, then $\overline{V}_G (\Gamma)
\ge \text{Vol} (\Gamma \setminus G / K) > 0$.
This implies a complete higher-dimensional
analogue of Gromov-Thurston Theorem [10]:

\proclaim{Theorem 2}  Let $M = \Gamma \setminus G /K$ be a
locally homogeneous manifold.
For any $N$, there are only finitely many $d \in \bbn$,
which are degrees of
maps $f: N \to M$.  In particular, any map $f: M \to M$
has degree 0 or $\pm 1$.
\endproclaim

In case $G = SL_2 (\bbc)$ this is the Gromov-Thurston theorem.
$G = \widetilde{SL(2, \bbr)}$
yields all Seifert fibrations, the
main ``complementary class'' to hyperbolic manifolds
according to Thurston's Geometrization Conjecture.

In section 3 we discuss the relation between
$V_{\widetilde{SL (2, \bbr)}} (\Gamma)$ and the generalized Hopf invariant.We
then reformulate this connection combinatirially, in spirit with
Quillen-Sullivan theory, which enables us to establish an effective estimate of
the Chern-Simons invariant in terms of
the number of simplices of a given triangulation of
$M$, a ``combinatorial complexity'' of $M$.
As a result, we get an explicit bound for the growth of topological complexity
for finite converings of $M$:

\proclaim{Theorem 3} Let $M$ be a plumbing of Seifert fibrations.  Then the
combinatorial complexity
of a $d$-fold covering $N$ of $M$ grows at least as $(\log d)^{1 - \e}, \ \e >
0$.
\endproclaim

I wish to thank a number of people for interesting discussion and useful
information:  Michel Boileau, Robert Brooks, Ilya Rips, Shicheng Wang,
Alexander
Lubotzky, Hyman Bass and Rostislav Grigorchuk.  The work on the paper was
initiated during my visits to R\"uhr-Universit\"at Bochum and
Universit\'e Paul
Sabatier, Toulouse.  I wish to  emphasize that the treatment of
Problem A is directly influenced by Gromov's work [9].

\demo{1. Axioms for volumes}

Let {\frak\$G} be the set of isomorphism classes of
finitely-presented groups.
A nonnegative function $\overline{V}: \text{\frak\$G} \to \bbr$ is
called an
upper volume if for $\Delta \subset \Gamma$ a subgroup of finite
index one has
$$ \overline{V} (\Delta) \le (\Gamma: \Delta) \overline{V} (\Gamma). $$
Similarly, a lower volume is a function yielding the opposite inequality.
A volume is a function for which an equality always holds.

Let $V: \text{\frak\$G} \to \bbr$ be a volume.  A group $\Gamma$ is
call distinctable
by $V$ if $V (\Gamma)$ is positive.  The importance of this
notion is seen form the
following proposition.

\proclaim{1.1. Proposition}  Let $V$ be a volume
function and let $\Gamma$ be
distinctable by $V$.  Then any two isomorphic
finite-index subgroups of $\Gamma$
have the same index.
\endproclaim

Let $V$ be a volume.  We say that $V$ is Hopfian, if for any
epimorphism $\Gamma \to \Pi$
one has $V (\Pi) \le V (\Gamma)$.

\proclaim{1.2. Proposition}  Let $V$ be a Hopfian volume,
and let $\Gamma$ be
distinctable by $V$. Then the image of any endomorphism $\vp:
\Gamma \to \Gamma$ is either all of $\Gamma$, or of
infinite index in $\Gamma$.  If $\Gamma$ is residually finite,
then $\vp$ is an isomorphism in the
first case.
\endproclaim

\demo{Proof}  Let $\Pi = \vp (\Gamma)$
and suppose $K (\Gamma: \Pi) < \infty$.
Then $V (\Gamma) \ge V (\Pi) = (\Gamma: \Pi) V (\Gamma) > V (\Gamma)$
a contradiction.
So either $\Pi = \Gamma$ or $(\Gamma: \Pi) = \infty$. Since any residually
finite group is Hopfian, the second statement also holds.

\demo{1.3}  Let $\underline{V}$ be a lower volume,
$\overline{W}$ is an upper volume and $V \le W$.  Define
$$ V (\Gamma) = \underset \Delta \rightarrow 1 \to {\underline{\lim}} \quad
\frac{\underline{V}(\Delta)}{[\Gamma: \Delta]} \quad,
\overline{W} (\Gamma) = \underset
\Delta \rightarrow 1 \to {\overline{\lim}} \quad
\frac{\overline{W} (\Delta)}{(\Gamma: \Delta)} $$
where ${\underline{\lim}}$, are taken over the lattice of finite index
subgroups.  We claim

\proclaim{Proposition 1.3}  $V$ and $W$ are volumes and $V \le W$.
\endproclaim

\demo{Proof}  One checks immediately that $0 \le V (\Gamma) \le W
(\Gamma) < \overline{W} (\Gamma) < \infty$.  Since for $\Pi \le
\Delta \le \Gamma,
\qquad \frac{\underline{V} (\Pi)}{[\Delta : \Pi]} = [\Gamma : \Delta] \
\frac{\underline{V} (\Pi)}{[\Gamma: \Pi)}, \quad V$ is volume and
similarly for $W$.

\demo{1.4 Example}  Define $V_\chi (\Pi) = \cases | \chi (\Gamma) |
&\text{if} \ \Gamma \ \text{is} \ F P \\
0 &\text{otherwise}. \endcases $  Then $V_\chi$ is a volume.

\demo{Proof}  Follows form [17], [7].

More interesting examples are given in the next section.

\subheading{2. Deficiency and rank volumes}  Let def$_+
(\Gamma) = \max (\text{def} (\Gamma), 1)$.

\demo{2.1}  Let $\underline{V}_d (\Gamma) = \text{def}_+ (\Gamma) - 1$
and $\overline{V}_r (\Gamma) = r (\Gamma) - 1$.  We claim that
$\underline{V}, \overline{W}$
yield the conditions of the Proposition 1.3.  Indeed, obviously
$\underline{V}_d (\Gamma) \le \overline{V}_r (\Gamma)$
since def$(\Gamma) \le r (\Gamma)$.  Moreover, for
$\Delta \le \Gamma$ of index
$d$ one has
$$ \align r (\Delta) - 1 &\le d (r (\Gamma) - 1) \\
\text{def} (\Delta) - 1 &\ge d (\text{def} (\Gamma) - 1) \endalign $$
[11].  So $\underline{V}_d$ is a lower volume and $\overline{V}_r$is an
upper volume.  By Proposition 1.3 we get volumes $V_d$ and $V_r$ satisfying
$$ \text{def}_+ (\Gamma) - 1 \le V_d (\Gamma) \le V_r (\Gamma) \le r (\Gamma) -
1 $$
In particular, if def$(\Gamma) \ge 2$, then $\Gamma$ is distinctable by $V_d$
and
$V_r$.  We are now in position to prove Theorem 1.

\demo{Proof of Theorem 1}  Since $\Gamma$ is distinctable by $V_d$, this
follows from Proposition 1.1.

\demo{2.2} $\underline{V}_r$ \underbar{is Hopfian}.  We claim that $V_r$ is a
Hopfian volume.  By definition $V_r (\Gamma) = \underset \Delta \rightarrow 0
\to {\overline{\lim}}
\frac{\overline{V}_r (\Delta)}{[\Gamma: \Delta]}$.  Let
$\Gamma \overset \vp \to
\twoheadrightarrow \Pi$ be onto and let $\Lambda \le \Pi$ be of index $d$.
Then put $\Delta = \vp^{-1} (\Lambda)$ so that
$[\Gamma: \Delta] = d$ and $\vp |_\Delta : \Delta
\to \Lambda$ is onto.  Since $r (\Lambda) \le r (\Delta)$ we get $V_r (\Gamma)
\ge V_r (\Pi)$.

\proclaim{Corollary}  Let $M$ be a compact manifold with residually finite
fundamental group of deficiency $\ge 2$.  Then any nonzero degree map $f: M \to
M$ induces
an automorphism of $\pi_1 (M)$.
\endproclaim

\demo{Proof}  Since deg$f \neq 0$, the image $f_\ast (\pi_1 (M))$ is of
finite index in $\pi_1 (M)$, by the standard argument.  Then Proposition 1.2
implies
that $f_\ast$ is an isomorphism.

\subheading{3. Representation and volumes}  Let $G$ be real semisimple Lie
group.  Let $K$ be the maximal compact subgroup of $G$ and let $X^n = G/ K$.
Fix a $G$-invariant volume form $\omega$ an $X$.  By Borel, see e.g. [14] one
has
an element Bor$(\omega) \in H^n_{\text{Cont}} (G)$.  If $\Gamma
\overset \rho \to
\rightarrow G$ is a representation one has a map
$B \Gamma \to B G^\delta$ and an induced map
$H^\ast (G^\delta, \bbr) \to H^\ast (\Gamma, \bbr)$.  Composing with a natural
map $H^n_{\text{Cont}} (G) \to H^n (G^\delta, \bbr)$ one gets an element
Bor $(\rho, \omega ) \in H^n (\Gamma, \bbr)$.  We refer to [14] for
direct geometric construction of this element and numerous applications.

Suppose we have a closed oriented manifold $M^n$ which is $K (\Gamma, 1)$.
Then
there is an element $[M] \in H_n (\Gamma, \bbz)$ and in fact $H_n (\Gamma,
\bbz) \approx \bbz$
is generated by it.  We put
$$ \text{Vol} (\rho) = (\text{Bor} (\rho, \omega), \ [M]) \in \bbr.$$
Now, we put
$$ V_G (\Gamma) = \max\limits_\rho | \text{Vol} (\rho) | $$
One needs an argument to show that the maximum is finite.  In fact,
the Cheeger-Simons rigidity [6] implies that Vol$(\cdot)$ is a locally
constant function on representation variety Hom$(\Gamma, G)$.  On
the other hand, the latter has finitely many components [9].  So $V_G (\Gamma)$
is finite.

Let $\Delta < \Gamma$ be of finite index and let $N \to M$ be the corresponding
covering.  The diagram
$$ \matrix N &\rightarrow &M \\
\| &  &\| \\
B \Delta &\rightarrow &\qquad \qquad B \Gamma \rightarrow B G^\delta \endmatrix
$$
shows that $V_G (\Delta) \ge (\Gamma : \Delta) V_G (\Gamma)$.  That is, $V_G
(\cdot)$ is a
lower volume.  Moreover, let $Q$ be any closed manifold and let $f: Q \to M$ be
a continuous map.  Then
$$ V_G (\pi_1 (Q)) \ge \text{deg} f \cdot V_G (\Gamma) $$
\demo{3.2}  Now, let $\Gamma$ be a torsion-free uniform lattice in $G$.  Then
[14]
$V_G (\Gamma) \ge \text{mes} ( \Gamma \setminus G / K)$, where the measure is
induced by $\omega$.  In particular,
$V_G (\Gamma)$ is positive.

Now we are ready to prove

\proclaim{Theorem 2}  Let $M = \Gamma \setminus G / K$ be a locally homogeneous
compact manifold.  For any compact oriented $Q, \ \dim Q = \dim M$, there
are only finitely many $d \in \bbz$ which are degrees of many $f: Q \to M$.
In particular, any selfmap $f: M \to M$ has degree 0 or $\pm 1$.
\endproclaim

\demo{Proof}  The inequality above implies
$$ \text{deg} f \le \frac{V_G (\pi_1 (Q))}{V_G (\Gamma)} $$
and since $V_G (\Gamma) > 0$, the RHS is finite.  The second statement follows
by application of the first statement to iterations of $f$, as usual.

\demo{Example} Let $G = S O (n, 1)$ so that $M$ is compact hyperbolic.  We
claim
that $V_G (\Gamma)$ is just a hyperbolic volume of $M$.  Indeed, for any
representation $\rho$ we have
$$ \text{Vol} (\rho) \le \text{Vol} (\rho_0) = \text{Vol} (M) $$
where $\rho_0$ is the natural representation ([14]).

In general, for any $M$ one has $V_G (\pi_1 (M)) \le \mu_n \| M \|$ where
$\| M \|$ is Gromov's simplicial volume and $\mu_n$ is the Milnor constant.
[14]

\demo{3.3.}  Now let $G = \widetilde{SL_2 (\bbr)}$ and let $M^3$ be a rational
homology sphere.  Let $\Gamma = \pi_1 (M)$ and let $\rho: \Gamma \to \widetilde
{SL_2(\bbr)}$ be
a representation.  We wish to interpret Vol $(\rho)$ as a
\underbar{generalized Hopf}\break \underbar{invariant}.
Recall that $\widetilde {SL_2 (\bbr)}$ acts on the hyperbolic plane
{\frak\$H}$^2$.
Let $V$ be the area form on {\frak\$H}$^2$.  Form a flat fiber bundle
$$ F = \tilde M \underset \pi_1 (M) \to \times \text{\frak\$H}^2$$
it carries a self-parallel fiber-like two-form $\nu$.

This extends naturally to a closed two-form an $F$, using the flat connection;
we
keep the notation $\nu$ for this form [14].  Now, for
any smooth section $s$ of
$F \to M$ put $\lambda = s^\ast \nu$. This is a closed two-form on $M$.  Since
$M$ is a rational homology sphere, $\lambda = d \k$ for some
$\k \in \Omega^1 (M)$.  Finally,
compute $H (s) = \int\limits_M \k \cdot \lambda$.  We claim that
independently on $s, H (X) = \text{Vol}
(\rho)$. Indeed, since {\frak\$H}$^2$ is contractible, all sections are
homotopic.  Now
the standard proof of the homotopy-invariance of Hopf invariant
(in Whitehead form, see\break [8])
applies in our situation and shows that $H (s)$ is independent
on $s$.  Next, consider the diagram
$$ \matrix \tilde M  \rightarrow &\widetilde{SL_2 (\bbr)} \\
&\downarrow \\
&\text{\frak\$H}^2 \endmatrix$$
where the horizontal map $\vp$ is any equivariant smooth map.  For
any left-invariant $\beta$ such that $d \beta = \pi^\ast \nu$ we have
$\beta \cdot
\pi^\ast \nu = \omega$.  Let $s$ be a section of ${\Cal F}$, coming form $\vp$,
then obviously $H (s) = \text{Vol}
(\rho)$.

Now, consider a triangulation of $M$ and let $\gamma$ be a simplicial
two-cochain, given by integration of $s^\ast \nu$ over simplices.  We
can compute $H (s)$ combinatorially, finding one-cochain $\alpha$ such
that $d \alpha = \gamma$ and then computing $(\alpha \cup \gamma, [ M])$.  Let
$s$ be a section, given by Thurston straightening
technique $([14] \ \ )$.
Then $(\gamma, \sigma) \le \pi$ for any two-simplex $\sigma$.

The map $d + \delta: C^1 (M, \bbr) \oplus \tilde C^3 (M, \bbr) \to
\tilde C^0 (M, \bbr) \oplus C^2 (M, \bbr)$ is an isomorphism since $M$ is
a rational homology sphere.  Let $c_i = \dim C^i$ be a number of $i$-simplices.
The matrix of $d + \delta$ is a $(c_1 + c_3 - 1) \times (c_1 + c_3 - 1)$
integer
invertible matrix, with all entries zero or one and the number of nonzero
entries in any
column or row does not exceed the
adjacency at the triangulation (the number of simplices of the next dimension,
adjacent to a given one), say $a$. So the sup-norm of $(d + \delta)^{-1}
(\gamma)$ is
bounded by $(\sqrt{a})^{c_1 + c_3 -1} (c_1 + c_3 -1) \cdot \pi$, by
Hadamard inequality.  It follows that
$$ | \text{Vol} (\rho) | = | (\alpha \cup \gamma, [M])
| \le ( \sqrt{a})^{c_1 +
c_3 - 1} (c_1 +
c_3 - 1) \cdot \pi^2 \cdot c_3. $$

Now, let $\e : N \to M$ be a $d$-sheeted covering.  For any triangulation of
$N$ we have
$$ d | \text{Vol} (\rho)| \le ( \sqrt{a'})^{c'_1 + c'_3 -1} (c'_1 + c'_3 -1)
\pi^2 c'_3, $$
or, assuming Vol$(\rho) > 0, \ \frac12 (c'_1 + c'_3 - 1)
\log a' \ge \log d + \text{const} (M) + O ( \log (c'_1 + c'_3) ).$
Now $c'_1 \le 2 c'_3$ and $a' \le c'_3$.  So
$$ \frac32 c'_3 \log c'_3 \ge \log d + \text{const} (m) + O (\log c'_3)). $$
Summing up, we have

\proclaim{Theorem 3}.  Let $M$ be a geometric non-hyperbolic homology sphere
(i.e. a
plumbing of Seifert fibrations).  Then the topological complexity
of its $d$-sheeted coverings grows at least as $(\log d)^{1 - \e}$, for
any $\e > 0$.
\endproclaim

\demo{Proof}  Any such $M$ admits a representation
$$ \rho : \pi_1 (M) \to \widetilde{SL_2 (\bbr)} $$
with Vol$(\rho) \neq 0$.

\subheading{4. Concluding remarks}

\demo{4.1. Volumes for pro-$p$ groups}  Let $G$ be a pro-$p$ group with
finitely generated homology groups. Then $\dim_{\bbf_p} (H^k (G, \bbf_p))$ is
a lower volume, that is, for $H < G$ an open subgroup.
$$ \dim_{\bbf_p} H^k (H, \bbf_p) \le [G: H] \dim_{\bbf_p} H^k (G, \bbf_p)). $$
Therefore $\mu_k (G) = \underset H < G \to {\overline{\lim}}
\frac{\dim H^k (H, \bbf_p)}{[G: H]}$ is
a finite volume.  If $\dim_{\bbf_p} H^1 (G) - \dim_{\bbf_p} H^2 (G) \ge 2$,
then $\mu_1 (G) > 0$.  We refer to [15] for details.
\pbf
\centerline{\bf References}

\item{[1]} H. Bass, Euler characteristic and characters of discrete groups.
Inv. Math., {\bf 28} (1976), 351--342.

\item{[2]} R. Brooks, W. Goldman, Volumes in Seifert space, Duke Math. J., {\bf
51} (1984), 529--545.

\item{[3]} R. Brooks, W. Goldman, The Godbillon-Wey invariant
of a transversely
homogeneous foliation, TAMS {\bf 286} (1984), 651--664.

\item{[4]} K. Brown, Euler characteristics of discrete groups and $G$-spaces,
Inv.
Math. {\bf 27} (1974), 229--264.

\item{[5]} J. Carlson, D. Toledo, Harmonic mappings of K\"ahler manifolds to
locally symmetric spaces, Publ. Math. IMES, {\bf 69} (1989), 173--201.

\item{[6]} J. Cheeger, J. Simons, Differential characters and geometric
invariants, in
Geometry and topology, J. Alexander and J. Karer, eds, LNM 1167, Springer,
1985, 50--80.

\item{[7]} I. Chiswell, Euler characterics of Groups, Math. Z., {\bf 147}
(1976), 1--11.

\item{[8]} E.Friedlender,P. Griffiths, J. Morgan, Homotopy theory and
differential forms, Mimeographed notes (1972).

\item{[9]} M. Gromov, Volume and bounded cohomology, Publ. Math.
IHES, {\bf 56} (1983), 5--99.

\item{[10]} M. Gromov, Hyperbolic manifolds according to
Thurston and Jorgensen, Sem. Bourbaki 1979/1980, 40--53.

\item{[11]} A. Lubotzky, Group presentation, $p$-adic analytic groups and
lattices in
$SL(2, \bbc)$, Annals of Mathematics, {\bf 118} (1983), 115--130.

\item{[12]} R. Lyndon, Problems in combinatorial group theory, Annals of Math.
Stud. {\bf 111} (1987), 3--33.

\item{[13]} A. Reznikov, Harmonic maps, hyperbolic cohomology and higher
Milnor inequalities, Topology, {\bf 32} (1993), 899--907.

\item{[14]} A. Reznikov, Rationality of secondary classes, J. Diff. Geom., to
appear.

\item{[15]} A. Reznikov, A spectral sequence for pro-$p$ group cohomology and
Adem inequalities, in preparation.

\item{[16]} J.Stallings,An extension theorem for Euler characteristics of
groups, preprint (1976).

\item{[17]}  C. T. C. Wall, Rational Euler characteristics, Proc. Cam. Phil.
Soc. {\bf 57} (1961), 182--183.

\item{[18]} S. Wang, Y. Q. Wu, Covering invariants and co-Hopficity of
3-manifold
groups, Proc. LMS, {\bf 68} (1994), 203--224.
\pmf
Institute of Mathematics
\pmf
The Hebrew University
\pmf
Givat Ram 91904, Jerusalem
\pmf
e-mail: simplex\@math.huji.ac.il

and

Max-Planck-Institut f\"ur Mathematik
\pmf
e-mail:reznikov@mpim-bonn.mpg.de
\bye